\def\draft{n}
\documentclass{amsart}
\usepackage[headings]{fullpage}
\usepackage{amssymb,epic,eepic,epsfig,amsbsy,amsmath}

\theoremstyle{plain}

\newtheorem{theorem}{Theorem}
\newtheorem{proposition}{Proposition}[section]
\newtheorem{lemma}[proposition]{Lemma}

\theoremstyle{definition}
\newtheorem{definition}[proposition]{Definition}

\newtheorem{question}{Question}

\theoremstyle{remark}

\newtheorem{remark}[proposition]{Remark}

\def\printname#1{
    \if\draft y
        \smash{\makebox[0pt]{\hspace{-0.5in}
            \raisebox{8pt}{\tt\tiny #1}}}
    \fi
}

\newcommand{\psdraw}[2]
         {\begin{array}{c} \hspace{-1.3mm}
    \raisebox{-4pt}{\epsfig{figure=draws/#1.eps,width=#2}}
    \hspace{-1.9mm}\end{array}}

\newlength{\standardunitlength}
\catcode`\@=11
\long\def\@makecaption#1#2{%
     \vskip 10pt

\setbox\@tempboxa\hbox{
       \small\sf{\bfcaptionfont #1. }\ignorespaces #2}%
     \ifdim \wd\@tempboxa >\captionwidth {%
         \rightskip=\@captionmargin\leftskip=\@captionmargin
         \unhbox\@tempboxa\par}%
       \else
         \hbox to\hsize{\hfil\box\@tempboxa\hfil}%
     \fi}
\font\bfcaptionfont=cmssbx10 scaled \magstephalf
\newdimen\@captionmargin\@captionmargin=2\parindent
\newdimen\captionwidth\captionwidth=\hsize
\catcode`\@=12

\newcommand{\qbinom}[2]{\text{$\left[\begin{array}{c}#1\\ #2\end{array}
\right]$}}

\def\lbl#1{\label{#1}\printname{#1}}

\def\BN{\mathbb N}
\def\BZ{\mathbb Z}

\def\BQ{\mathbb Q}
\def\BR{\mathbb R}
\def\BC{\mathbb C}

\def\a{\alpha}

\def\La{\Lambda}
\def\l{\lambda}

\def\e{\epsilon}

\def\d{\delta}
\def\b{\beta}

\def\mat#1#2#3#4{\left(
\begin{matrix}
 #1 & #2  \\
 #3 & #4
\end{matrix}
\right)}

\def\longto{\longrightarrow}

\def\bk{\mathbf{k}}

\def\k{\kappa}
\def\ev{\mathrm{ev}}
\def\za{z_{\alpha}}
\def\zb{z_{\beta}}
\def\zk{z_{\kappa}}
\def\pt{\partial}
\def\GVC{\rm{GVC}}
\def\VC{\rm{VC}}
\def\fg{\mathfrak{g}}
\def\SL{\mathrm{SL}}

\begin{document}


\title[On the volume conjecture for small angles]{
On the volume conjecture for small angles}

\author{Stavros Garoufalidis}
\address{School of Mathematics \\
         Georgia Institute of Technology \\
         Atlanta, GA 30332-0160, USA \\
         {\tt http://www.math.gatech} \newline {\tt .edu/$\sim$stavros } }
\email{stavros@math.gatech.edu}
\author{Thang TQ L\^e}
\address{School of Mathematics \\
         Georgia Institute of Technology \\
         Atlanta, GA 30332-0160, USA}
\email{letu@math.gatech.edu}

\thanks{The authors were supported in part by National Science Foundation. \\
\newline
1991 {\em Mathematics Classification.} Primary 57N10. Secondary 57M25.
}

\date{February 4, 2005 \hspace{0.5cm} First edition: February 4, 2005.}

\dedicatory{Dedicated to Louis Kauffman on the occasion of his 60th birthday}

\begin{abstract}
Given a knot in 3-space, one can associate a sequence of Laurrent polynomials,
whose $n$th term is the $n$th colored Jones polynomial. The Generalized 
Volume Conjecture states that the value of the $n$-th colored Jones
polynomial at $\exp(2 \pi i \a/n)$ is a sequence of complex numbers that
grows exponentially, for a fixed real angle $\a$. 
Moreover the exponential growth rate 
of this sequence is proportional to the volume of the 3-manifold obtained
by $(1/\a,0)$ Dehn filling. In this paper we will prove that (a) for every
knot, the limsup in the hyperbolic volume conjecture is finite and bounded
above by an exponential function that depends on the number of crossings.
(b) Moreover, for every knot $K$ there exists a positive real number $\a(K)$ 
(which depends on the number of crossings of the knot)
such that the Generalized Volume Conjecture holds for $\a \in [0,
\a(K))$. Finally, we point out that a theorem of Agol-Storm-W.Thurston
proves that the bounds in (a) are optimal, given
by knots obtained by closing large chunks of the weave.
\end{abstract}

\maketitle

\tableofcontents


\section{Introduction}
\lbl{sec.intro}

\subsection{The volume conjecture}
\lbl{sub.volume}

The {\em Volume Conjecture} (in short, \VC ) 
connects two very different approaches to knot theory,  
namely Topological Quantum Field Theory and Riemannian
(mostly Hyperbolic) Geometry. 

In its fundamental form, the \VC\ states that for every hyperbolic knot $K$
in $S^3$ we have:
\begin{equation}
\lbl{eq.VC}
\lim_{n \to \infty} \frac{\log|\ev_{1,n}(J_K(n))|}{n}=
\frac{1}{2 \pi}\, V_K,
\end{equation}
where 
\begin{itemize}
\item
$\ev_{\a,n}(f)$ denotes the evaluation of a rational function $f(q)$ at
$q=e^{2 \pi i \a/n}$, 
\item
$J_K(n) \in \BZ[q^{\pm}]$ is the {\em Jones polynomial} of a knot 
{\em colored} with the 
$n$-dimensional irreducible representation of $\mathfrak{sl}_2$,
normalized so that it equals to $1$ for the unknot (see \cite{J, Tu}), and 
\item 
$V_K$ is the {\em volume} of the knot complement $S^3-K$, using the unique
complete hyperbolic metric; \cite{Th}.
\end{itemize}
The \VC\ 
was formulated  in this form by 
H. and J. Murakami \cite{MM}, who reinterpreted an earlier version 
due to Kashaev, \cite{K}.

It is natural to ask what happens when we use evaluations $\ev_{\a,n}$
of the colored Jones polynomial for other complex numbers $\a$. In \cite{Gu},
Gukov proposed a {\em Generalized Volume Conjecture} (in short, \GVC ),
which states that for every hyperbolic knot $K$ in $S^3$ 
and every irrational $\a$ near $1$ (or $\a=1$), we have: 
\begin{equation}
\lbl{eq.GVC}
\lim_{n \to \infty} \frac{\log|\ev_{\a,n}(J_K(n))|}{n}=
\frac{1}{2 \pi}\, V_K(1/\a),
\end{equation}
where
\begin{itemize}
\item
$V_K(1/\a)$ is the {\em volume} of the $(1/\a,0)$ Dehn filling of the 
knot complement $S^3-K$ (that is, the Dehn filling corresponding to 
$M^{2 \pi i \a} L^0=1$ where $(M,L)$ is the meridian and longitude of the 
knot); see \cite{Th}. 
\end{itemize}
The \GVC\
is mostly about hyperbolic knots and hyperbolic Dehn fillings. In case a 
manifold (with or without boundary) is not hyperbolic, we will declare its
corresponding Gromov-Thurston volume to be zero.

There are two rather independent parts in the \GVC :
\begin{itemize}
\item[(a)]
To show that the limit exists in \eqref{eq.GVC},
\item[(b)]
To identify the limit with the volume of the corresponding Dehn filling.
\end{itemize}

At the moment, the GVC is known for the $4_1$
knot and certain values of $\a$; see Murakami, \cite{M1}.

In the following, we will refer to the parameter $\a$ in the \GVC\ as 
{\em the angle}, making contact with standard terminology from
hyperbolic geometry.

One may further ask what happens to the \GVC\ when the angle $\a$ is small.
For example, when $\a=0$, then $\ev_{0,n}(J_K(n))=1$ for all knots $K$ and
integers $n$, thus the corresponding limit on the left hand side of
\eqref{eq.GVC} vanishes. To attach a meaning to the right hand side of
\eqref{eq.GVC}, we need to use the volume of an appropriate $\SL_2(\BC)$
representation of the knot complement. Unfortunately, for small angles
$\a$, the corresponding Dehn fillings are not hyperbolic. But fortunately,
there is a natural $\mathrm{SL}_2(\BC)$ representation to consider when
$\a=0$, namely the {\em trivial} one. This is in agreement with physics,
where the case $\a=0$ is a classical limit of a quantum theory, and corresponds
to the only $\SL_2(\BC)$ flat connection on $S^3$, namely the trivial one.
When $\a$ is small and real, we will {\em define} $V_K(1/\a,0)$ to be the 
volume of the {\em reducible} $\SL_2(\BC)$ representation 
$$ 
\rho_{\a}: \pi_1(S^3-K) \longto \SL_2(\BC),
\qquad \rho_{\a}(\text{M})= \mat {e^{2 \pi i \a/n}} 0 0
{e^{-2 \pi i \a/n}}. 
$$
(where $M$ is a meridian). A simple calculation shows that $V_K(1/\a,0)=0$
for $\a$ small and positive.

Now, we can formulate the \GVC\ for small positive angles $\a$. 

Thus, we have formulated a \GVC\ for $\a$ near $0$ and $\a$ near $1$, using
representations near the trivial one and a discrete faithful, respectively.
How can we connect and explain our choices for other angles $\a$? A natural
answer to this question requires analyzing asymptotics of solutions of 
difference equations with a parameter. This is a different subject that
we will not discuss here; instead we will refer the curious reader to
\cite{GG}, and forthcoming work of the first author.

\subsection{Upper bounds and confirmation of the volume 
conjecture for small angles}
\lbl{sub.results}

Our results are the following:

\begin{theorem}
\lbl{thm.1}
For every knot $K$ with $c+2$ crossings and every $\a >0$, we have
$$
\limsup_{n \to \infty} \frac{\log|\ev_{\a,n}(J_K(n))|}{n}
\leq c \log 4.
$$
\end{theorem}

\begin{theorem}
\lbl{thm.2}
For every knot $K$, 
there exists a positive angle $\a(K)>0$ such that
the \GVC\ holds 
for all $\a \in [0, \a(K))$.
\end{theorem}

In fact, the proof of Theorem \ref{thm.2} reveals that we can take $\a(K)$
to be a function that depends on the number of crossings of $K$.

Notice that the Hyperbolic Volume Conjecture is the case of $\a=1$, which
corresponds to a complete hyperbolic structure. On the other hand, Theorem
\ref{thm.2} deals with small cone-angle fillings, which are expected to be
{\em spherical structures} of zero volume.

Theorem \ref{thm.1} follows easily from a stronger result. Observe that
when $|q|=1$, then
$$
|J_K(n)(q)| \leq ||J_K(n)||_1
$$
where the $l^1$-{\em norm} of a Laurrent polynomial $f=\sum_k c_k 
q^k$ is given by:
$$
||\sum_k c_k q^k||_1 = \sum_k |c_k|.
$$

Then, we have the following:

\begin{theorem}
\lbl{thm.3}
For every knot $K$ of $c+2$ crossings and every $n$ we have:
\begin{equation}
\lbl{eq.l1}
||J_K(n)||_1 \leq n^c 4^{c n}.
\end{equation}
\end{theorem}

\subsection{Relation with hyperbolic geometry, and optimal bounds}
\lbl{sub.hypgeom}

When $\a=1$, the upper bound in Theorem \ref{thm.1} is not optimal,
and does not reveal any relationship between the $\limsup$ and hyperbolic
geometry. Our next theorem fills this gap.

\begin{theorem}
\lbl{thm.4}
For every knot $K$ with $c+2$ crossings we have
$$
2 \pi \limsup_{n \to \infty} \frac{\log|\ev_n (J_K(n))|}{n}
\leq  v_8 c,
$$
where $\ev_n=\ev_{1,n}$, 
$v_8=8 \La(\pi/4) \approx 3.66386$ is the volume of the regular
ideal octahedron.
\end{theorem}

Using an ideal decomposition of a knot complement by placing one octahedron
per crossing, it follows that for every knot $K$ with $c+2$ crossings, we have
\begin{equation}
\lbl{eq.vol8}
V(K) \leq v_8 c.
\end{equation}
On the other hand, if the volume conjecture holds for $\a=1$, then
$$
2 \pi \limsup_{n \to \infty} \frac{\log|\ev_n (J_K(n))|}{n}
= V(K) \leq v_8 c.
$$
One may ask whether \eqref{eq.vol8} (and therefore, whether the
bound in Theorem \ref{thm.4}) is optimal. Optimality is at first sight
surprising, since it involves all knots (and not just alternating ones)
and their number of crossings (which carries little known geometric 
information).
In conversations with I.Agol and D.Thurston, it was communicated to us
that the upper bound in \eqref{eq.vol8} is indeed optimal. Moreover a class 
of knots that achieves (in the limit) the optimal ratio of volume 
by number of crossings is obtained by
taking a large chunk of the following {\em weave}, and closing it up to a knot:
$$
\psdraw{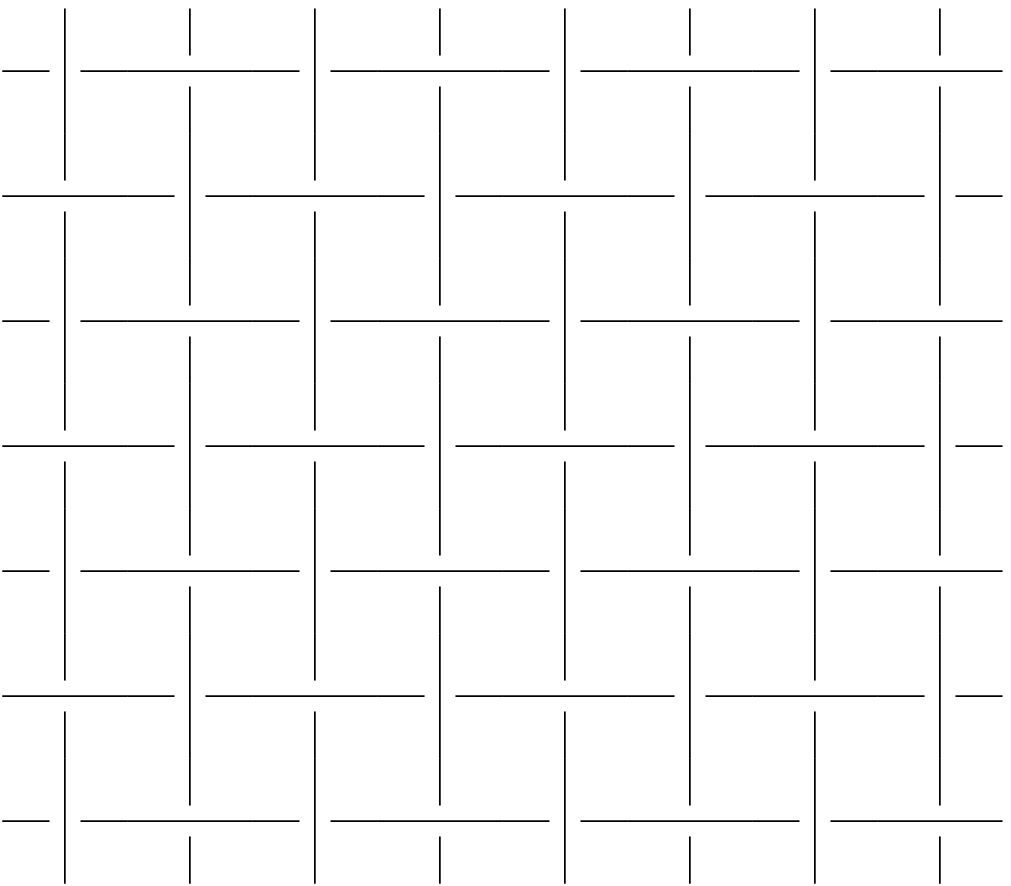}{2in}
$$
The complement of the weave has a complete hyperbolic structure associated
with the {\em square tessellation} of the Euclidean plane:
$$
\psdraw{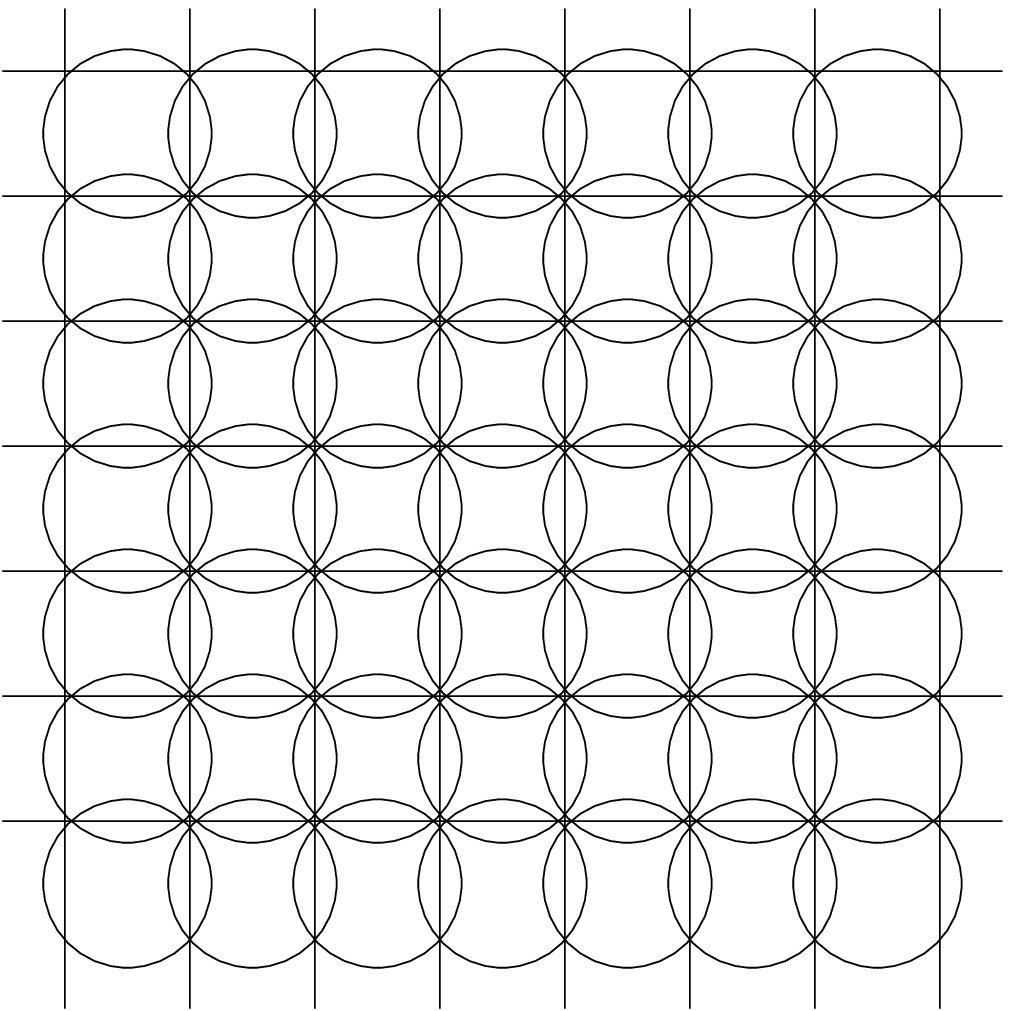}{2in}
$$
Optimality follows along similar lines as the Appendix of \cite{La},
using a stronger estimate for the lower bound of the volume of Haken
manifolds, cut along an incompressible surface. The stronger statement
is the following result which will appear in subsequent work of 
Agol-Storm-W.Thurston, \cite{AST}. Its proof uses, among other things,
work of Perelman.

\begin{theorem}
(\cite{AST})
If M is a hyperbolic finite volume 3-manifold containing a properly 
imbedded orientable, boundary incompressible, incompressible surface $S$, 
then
$$
V(M) >= V(\mathrm{Guts}(M-\mathrm{int(nbd(}S))),
$$
where $V$ stands for volume, and 
the $\mathrm{Guts}$ terminology are defined in \cite{Ag}.
\end{theorem}

The reader may compare \eqref{eq.vol8} 
with the following result of Agol-Lackenby-D.Thurston \cite{La}:

\begin{theorem}
\lbl{thm.ALT}
If $K$ is an alternating knot with a planar projection with $t$ twists, then
$$
v_3 (t-1)/2 <  V(K) < 10 v_3 (t-1),
$$
where $v_3=2 \La(\pi/3)\approx 1.01494$ is the volume of the regular ideal
tetrahedron. 
Moreover, the class of knots obtained by Dehn filling on the {\em chain link}
has asymptotic ratio of volume by twist number equal to $10 v_3$.
The corresponding tessellation of the Euclidean plane is given by the
{\em star of David}.
\end{theorem}

\subsection{The main ideas}
\lbl{sub.ideas}

Our results have quick proofs, and require only a small dose of
elementary analysis, and an appropriate view of these powerful quantum
invariants of knots.

We already saw how Theorem \ref{thm.1} follows from the stronger Theorem
\ref{thm.3}.

To prove Theorem \ref{thm.3}, we will make use of a state-sum definition
of the colored Jones polynomials, where the local weights are given by
$R$-matrices. In \cite{GL}, we used the basic fact that the local weights
are $q$-{\em hypergeometric}, 
in order to deduce that the sequence of colored Jones
polynomials satisfies a {\em linear recursion relation} (the recursion
depends on the knot, of course).
In our case, we will focus on the fact that the local weights 
take values in $\BZ[q^{\pm}]$
in order to give elementary estimates for their $l^1$-norm. 

Theorem \ref{thm.2} is trickier. Among other things, the proof uses a key
integrality property (due to Habiro)
of the cyclotomic transform of the colored Jones function. To the best of our
knowledge, this is a first application of the cyclotomic transform.
In more detail, 
we will separate out the dependence
of the color in the colored Jones polynomial $J_K(n)$, 
via the cyclotomic transform. This replaces $J_K(n)$ by a sequence $C_K(n)$ 
(for $n \in \BN$) of rational functions. 
Habiro proved that $C_K(n)$ are actually polynomials. 
Using the $l^1$-estimates on $J_K(n)$, one can deduce only a weak estimate
for the Mahler measure of $C_K(n)$. Due to the specific structure of the
inverse cyclotomic transform, one can prove a useful $l^1$-estimate
for $C_K(n)$. Using that estimate, and some elementary analysis,
it is easy to finish the proof of Theorem \ref{thm.2}.

In Section \ref{sec.qholonomic}, we exploit the fact that
the colored Jones function and the cyclotomic function
is a solution to a linear $q$-difference equation. Using elementary methods,
we give quadratic degree bounds for solutions of $q$-difference equations,
and exponential bounds for the $l^1$ norms of solutions of integral
(in the sense of Laurrent polynomials)  
$q$-difference equations. Based on experimental evidence,
we conjecture that the cyclotomic function is a solution of an integral
$q$-difference equation.

The upper bounds for the limsup in Theorem \ref{thm.1} are not sharp for 
$\a=1$, since they are obtained by $l^1$ estimates of the local weights.
Using the fact that the local weights are given by a ratio of 
$5$ quantum factorials, and the fact that the asymptotics of quantum factorials
are governed by the Lobachevsky function, in Section \ref{sec.hypgeom}
we give better bounds, stated in Theorem \ref{thm.4}, which are linear in the
number of crossings, and involve the volume of the regular ideal octahedron.
In addition, we conjecture that our improved bounds are optimal.

Finally, in two appendices we discuss the Volume Conjecture for 
the Borromean rings, and the Generalized Volume Conjecture for torus knots.

\section{Proof of Theorem \ref{thm.3}}
\lbl{sec.thm3}

As we mentioned above, we will make use of an $R$-matrix 
state sum definition of the colored Jones polynomial $J_K(n)$,  
discussed, for example, in \cite[Sec.3]{GL}. 
Consider a long knot $K'$ whose closure
is $K$, and fix a positive integer $n$. To compute $J_K(n)$, we follow
the following algorithm:

\begin{itemize}
\item
Assign angle variables $k$ at each crossing, and let $\bk=(k_1,k_2, \dots)$.
\item
Color each part-arc of the knot projection such that around each crossing
the color of the part-arcs is given by:
$$
\psdraw{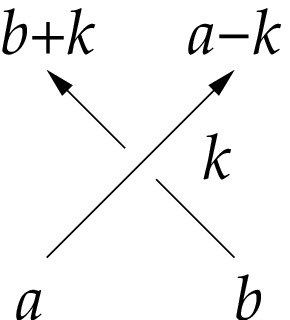}{0.5in} 
\qquad
\psdraw{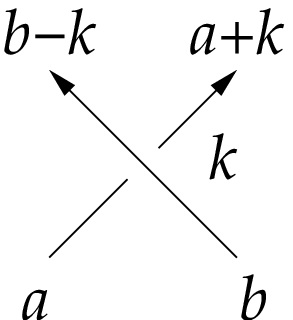}{0.5in} 
$$
There is a unique coloring of the part-arcs such that the two broken
part-arcs have color $0$.
\item
Assign local weights $f(n;a,b,k) \in \BZ[v^{\pm 1/2}]$ at each crossing,
and for a fixed coloring $\bk$, let $F(n,\bk)$ denote the product over
all crossings of the corresponding local weights, times $v$ raised to 
a linear form on $n,\bk$.
\item
Form the sum
\begin{equation}
\lbl{eq.statesum}
J_K(n)=\sum_{\bk} F(n,\bk).
\end{equation}
\end{itemize}
Strictly speaking, in \cite{GL} we discussed the above algorithm for knots 
which are closures of braids, but a similar algorithm works for planar
projections of knots as well. This follows from the following figure:
$$
\psdraw{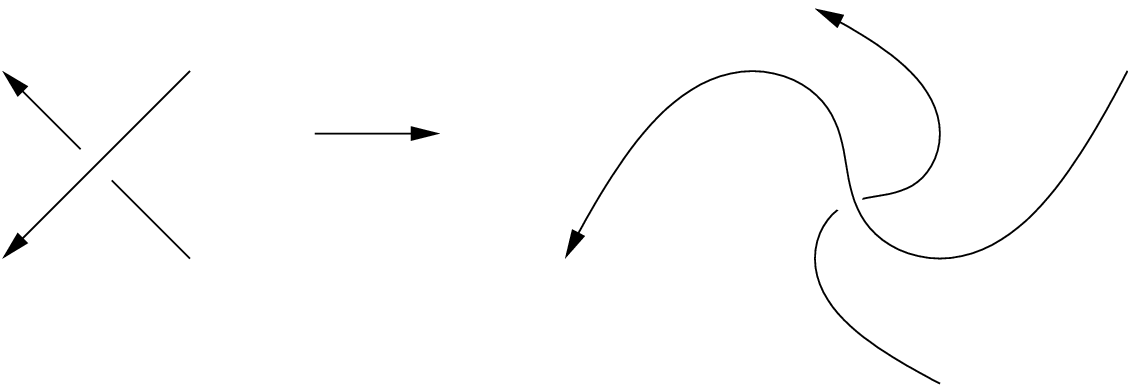}{2in}
$$
which moves crossings (positive or negative) into standard upright position,
by an isotopy that creates local minima/maxima and no further crossings.
The local minima/maxima give rise to an additional multiplicative factor 
(a monomial in $v$ raised to a linear form in $n,\bk$)  
in $F(n, \bk)$, and does not affect the estimates below.

In \cite{GL} we used the fact that the local weights $R_\pm:\BZ^5 \to 
\BZ[v^{\pm 1/2}]$, are $q$-holonomic functions, in order to deduce from
first principles that $J_K$ is $q$-holonomic, and thus satisfies a linear
recursion relation.

In our case, we will make use of the specific integral form of the local
weights in order to deduce our result. Let us recall the specific form of the
local weights, from \cite[Sec.3]{GL}.
\begin{eqnarray*}
\lbl{eq.R+-}
R_+(n;a,b,k) &:=&  (-1)^k v^{-((n-1-2a)(n-1-2b)+k(k-1))/2}
\qbinom{b+k}{k} \{n-1+k-a\}_k, \\
R_- (n;a,b,k) &:=&  v^{((n-1-2a-2k)(n-1-2b+2k)+k(k-1))/2}
\qbinom{a+k}{k} \{n-1+k-b\}_k,
\end{eqnarray*}
where $v=q^{1/2}$, and 
for $a,b \in \BN$, we define the $q$-{\em integers}, $q$-{\em factorial}
and $q$-{\em binomial} coefficients by:
\begin{equation}
\lbl{eq.qfactorial}
\{a\}:= v^a -v^{-a},
 \qquad \{a\}! := \prod_{i=1}^{a} \{a\},
\qquad
\{a\}_b := \frac{\{ a \}!}{\{ a-b \}!},
\qquad
\qbinom{a}{b} := \frac{\{ a \}!}{\{ b \}! \{ a-b \}!}
\end{equation}

For a fixed coloring $\bk$, the colors at the part-arcs of $K'$ are linear
forms on $\bk$ and $n$, and the contribution is nonzero only when all these 
linear forms are nonnegative and less than $n$. In particular, each of the
angle variables $k$ has to be $0 \leq k < n$.  
Observe that the $l^1$ norm satisfies the inequalities
\begin{equation}
\lbl{eq.ineq}
||f+g||_1 \leq ||f||_1+||g||_1, \qquad ||fg||_1 \leq ||f||_1 ||g||_1,
\end{equation}
Moreover, 
$\qbinom{m}{k} \in \BN[v^{1/2}]$ and 
$$
\qbinom{m}{k}=\binom{m}{k} \leq 2^m \leq 2^n
$$
if $m \leq n$. Moreover,
$$
||\{n-1+k-a\}_k||_1 \leq 2^k \leq 2^n
$$
for $k \leq n$. Combining with the above formulas for $R_{\pm}$, it follows
that $||R_{\pm}(n;a,b,k)||_1 \leq 4^n$. 
Since
$$
R_{\pm}(n;0,0,k)=\d_{k,0},
$$ 
and since we can always choose a breaking of a knot so that the broken arc
ends under the first crossing and out of the last crossing, it follows that
we can ignore at least two angle variables corresponding to
the broken part-arcs. Thus, 
$||F(n,\bk)||_1 \leq 4^{c n}$, where $c+2$ is the number of crossings
of $K$. Since $k$ has to be $0 \leq k <n$, \eqref{eq.statesum} implies that
we have at most $n^c$ choices for the angle variables. The result follows.
\qed

\begin{remark}
\lbl{rem.3}
There are several other formulations of the $n$-th colored Jones polynomial,
for example coming from the Kauffman bracket skein module. Unfortunately,
using the Kauffman bracket skein module formulation, it is hard to prove
Theorem \ref{thm.3}, since the number of crossings of an $n$-parallel of
a knot with $c+2$ crossings is $(c+2)n^2$, a quadratic function of $n$.
Nevertheless, the Kauffman bracket skein module can give good estimates
of the min and max degree (and their difference, the {\em span}) 
of the $n$-th colored Jones polynomial. 
Compare with \cite[Prop.1.2]{Le} of the second author, who observes that
$$
\text{span}(J_K(n)) \leq c n^2 + O(n).
$$
\end{remark}

\section{Proof of Theorem \ref{thm.2}}
\lbl{sec.thm2}

\subsection{A reduction of Theorem \ref{thm.2}}
\lbl{sub.reduction}

The proof of Theorem \ref{thm.2} will use the cyclotomic expansion of the 
colored Jones polynomial, introduced by Habiro, and an improved Mahler-type 
estimated, communicated to us by D. Boyd.

\begin{definition}
\lbl{def.cyclo}
Given a function $f: \BN \longto \BQ(q)$, we define
its {\em cyclotomic transform} $Cf: \BN \longto \BQ(q)$ by:
\begin{equation}
\lbl{eq.cyclo}
Cf(n)=\sum_{k=0}^\infty C(n,k) f(k)
\end{equation}
where $C(n,k)$ is given by:
\begin{eqnarray*}
\lbl{eq.cyclokernel}
C(n,k) &:=& \frac{1}{q^{n/2}-q^{-n/2}}
\prod_{j=n-k}^{n+k} (q^{j/2}-q^{-j/2}) \\ 
& = & \prod_{j=1}^k (( q^{n/2}-q^{-n/2})^2 - (q^{j/2}-q^{-j/2})^2) \\
& = & \prod_{j=1}^k (( q^{n/2}+q^{-n/2})^2 - (q^{j/2}+q^{-j/2})^2). 
\end{eqnarray*}
The cyclotomic transform has an inverse $C^{-1}f:\BN \longto \BQ(q)$
defined by:
\begin{equation}
\lbl{eq.cycloinverse}
C^{-1}f(n)=\sum_{k=0}^\infty R(n,k) f(k)
\end{equation}
where $R(n,k)$ is given by:
\begin{equation}
\lbl{eq.cyclokernelinverse}
R(n,k)= (-1)^{n-k} \frac{\{2k\} }{\{2n+1\}![2n]
}\qbinom{2n}{n-k},
\end{equation}
where for $a \in \BN$, we define:
$$
[a]:= \frac{\{a\}}{\{1\}}=\frac{v^{a/2}-v^{-a/2}}{v^{1/2}-v^{-1/2}}
$$
(and $v=q^{1/2}$). 
Notice that $R(n,k)=C(n,k)=0$ for $k>n$, so the above sums are finite.
\end{definition}

\begin{definition}
\lbl{def.cyclofunction}
If $J_K:\BN \longto \BZ[q^{\pm}]$ denotes the colored Jones function of
a knot $K$, we define the {\em cyclotomic function} of $K$ by $C_K=C J_K$.
\end{definition} 

A key result of Habiro is that the cyclotomic function of a knot takes values
in $\BZ[q^{\pm}]$; see \cite{H1}.

\begin{proof}(of Theorem \ref{thm.2})
Let us apply the cyclotomic transform to $J_K(n)$, in order to isolate
the dependence of the color. By the above definition, we have:
\begin{equation}
\lbl{eq.J2C}
J_K(n)=\sum_{k=0}^n C(n,k) C_K(k)=1+ \sum_{k=1}^n C(n,k) C_K(k). 
\end{equation}
where
\begin{equation}
\lbl{eq.C2J}
\{2n+1\}![2n] C_K(n)=\sum_{k=0}^n (-1)^{n-k} \{2k\} 
\qbinom{2n}{n-k} J_K(k)
\end{equation}

Let us assume for the moment the following:

\begin{theorem}
\lbl{thm.boundC}
For every knot $K$ we have:
\begin{equation}
\lbl{eq.LC}
||C_K(n)||_1 \leq e^{C n + O(\log n)}
\end{equation}
\end{theorem}

Here, and below, the $O(f(n))$ notation means that the error is bounded
by a constant times $f(n)$.

The reader may wonder what we gained starting from $J_K$, going to $C_K$,
and then back to $J_K$. The point of Equation \eqref{eq.J2C} is that
it separates the dependence of $J_K(n)$ on the color $n$, and Theorem
\ref{thm.boundC} gives exponential bounds for the $l^1$ norm of 
the polynomials $C_K(n)$.

For $q=e^{2 \pi i \a/n}$ and $\a$ small and positive,
and $0 < k < n$, $C(n,k)$ becomes small. Using Lemma 
\ref{lem.estimate} below and Theorem \ref{thm.3}, it follows that
when 
$0 < k < n$, then
\begin{eqnarray*}
|\ev_{\a,n}(C_K(k) C(n,k))| 
& \leq &  e^{C k + O(\log k)} |3 \sin(\pi \a)|^{2k} \\
& = & e^{ C'(\a) k + O(\log k)},
\end{eqnarray*}
where $C'(\a):=C + 2 \log(3 |\sin(\pi \a)|)$.
Now, choose $\a$ small enough 
so that $C'(\a)<0$,
and then choose $k_0=k_0(\a)$ so that 
$C'(\a) k + O(\log k)<
C'(\a)/2 \, k$ for $k \geq k_0$.
It follows that for $k_0 < k < n$, we have
$$
|\ev_{\a,n}(C_K(k) C(n,k))| \leq e^{C'(\a)/2 \, k}
$$
and the last term is in absolute value less than $1$.
Moreover, for $0 < k < k_0$, we have 
$$
\lim_{n \to \infty} \ev_{\a,n}(C_K(k) C(n,k))=0.
$$
This and Equation \eqref{eq.J2C} implies Theorem \ref{thm.2}. 
\end{proof}

\begin{lemma}
\lbl{lem.estimate}
If $0 < \a < \pi/6$, 
then for every $0 \leq k < n$
we have:
$$
|\ev_{\a,n}(C(n,k))| \leq |3 \sin(\pi \a)|^{2k}.
$$
\end{lemma}

\begin{proof}
Evaluating at $q=e^{2 \pi i \a/n}$, we have:
$$
|\ev_{\a,n}(C(n,k))|=\prod_{k=1}^n|q^n-q^k||q^n-q^{-k}|.
$$
Since $0 < \a < \pi/6$ and $q^n=e^{2 \pi i \a}$ it follows that for all 
$0 < k < n$ we have: 
$$
|q^n-q^k| \leq |q^n-q| < |q^n-1|=|2 \sin(\pi \a)|
$$
and 
$$
|q^n-q^{-k}| \leq |q^n-q^{-1}| < |q^{n+1}-1|=|2 \sin(\pi \a(n+1)/n)|
\leq |3 \sin(\pi \a)|.
$$
The result follows.
\end{proof}

It remains to prove Theorem \ref{thm.boundC}.

\subsection{Proof of Theorem \ref{thm.boundC}}
\lbl{sub.thmboundC}

Consider Equation \eqref{eq.C2J}. First of all it gives:
\begin{equation}
\lbl{eq.degC}
\text{span} C_K(n) \leq C''' n^2.
\end{equation}
Since $||\qbinom{2n}{n-k}||_1=\binom{2n}{n-k} \leq 2^{2n} $, Equation 
\eqref{eq.C2J} and Theorem \ref{thm.3} imply that
\begin{equation}
\lbl{eq.L1C}
||\{2n+1\}![2n] C_K(n)||_1 \leq n^c 4^{(c+1)n}.
\end{equation}
A priori, this estimate is weak and implies an exponential
upper bound on the Mahler measure of $C_K(n)$, and a doubly exponential
upper bound on the $l^1$ norm of $C_K(n)$. Let us digress a bit and
discuss this in detail.

Recall that the $l^2$-norm
$$
||f||_2:=(\sum_k |a_k|^2)^{1/2}
=\left(\int_0^1 |f(e^{2 \pi i t})|^2 dt \right)^{1/2}
$$   qf a polynomial
$f=\sum_k a_k q^k$ also satisfies the inequalities of \eqref{eq.ineq}.
However, neither the $l^1$ nor the $l^2$ norm are multiplicative.
Mahler introduced a {\em measure}
$$
M(f)=\exp \int_0^1 \log|f(e^{2 \pi i t})| dt
$$
which although it is not a norm, it is by definition multiplicative:
$$
M(fg)=M(f)M(g).
$$
The next proposition summarizes how the Mahler measure compares with the
$l^1$ and $l^2$ norms:

\begin{proposition}
\lbl{prop.mahler}
If $f$ is a Laurrent polynomial of degree $d$, then
\begin{eqnarray*}
M(f) & \leq & ||f||_2 \leq ||f||_1 \\
||f||_1 & \leq & 2^d M(f) \\
M(f) & = & 1 \quad \text{if $f$ is a product of cyclotomic polynomials}.
\end{eqnarray*}
\end{proposition} 

Returning to the proof of Theorem \ref{thm.2}, Equation \eqref{eq.L1C}
and Proposition \ref{prop.mahler} imply that
$$
||C_K(n)||_1 \leq n^c 4^{(c+1)n} 2^{\text{deg}C_K(n)}.
$$
We will see later that $\text{deg} C_K(n)=O(n^2)$. Thus, the above estimate
is exponential in $n^2$.

At this point, we will use the following theorem, communicated
to us by D. Boyd.

\begin{theorem}
\lbl{thm.boyd}(Boyd)
If $f(q)$ is a polynomial that satisfies 
\begin{itemize}
\item
$\text{deg}(f(q))= C n^2 + O(n)$ and 
\item
$||(1-q)(1-q^2)...(1-q^n) f(q)||_1 \leq e^{C' n + O(\log n)}
$
then
\end{itemize}
$$
||f||_1 \leq e^{C'' n + O(\log n)}. 
$$
\end{theorem}

\begin{proof}
Consider the polynomial
$
g(q)=f(q)(1-q)\dots(1-q^n)
$
of degree at most $Cn^2 + n(n+1)/2+O(n)$.
Thus, we may write
$$
g(q)=\sum_{k=0}^{(C+1/2) n^2 + O(n)} a_k q^k.
$$
Let us expand
$$
\frac{1}{\prod_{k=1}^n (1-q^k)}=\sum_{k=0}^\infty c_k q^k
$$
where $c_k \in \BN$. We can obtain an upper bound for the growth
of $c_k$ as follows.
Consider
$$
\frac{1}{\prod_{k=1}^\infty (1-q^k)}=\sum_{k=0}^\infty p_k q^k
$$
where $p_n$ is the {\em number of partitions} of $n$.
Using growth rate of $p_n$ (see \cite{An}), it follows that
\begin{equation}
\lbl{eq.partitions}
0 \leq c_k \leq p_k = e^{ \pi \sqrt{2/3} 
\sqrt{k}+O(1)}.
\end{equation}
The important thing is that the growth rate of $p_k$ involves $\sqrt{k}$.
Moreover, we have:
\begin{eqnarray*}
f(q) &=& \frac{g(q)}{(1-q)\dots(1-q^n)} \\
&=& \sum_{k=0}^{C n^2+O(n)} d_k q^k
\end{eqnarray*}
where
$$
d_k=\sum_{i=0}^k a_i c_{k-i}.
$$
Since $|a_i| \leq ||g(q)||_1 \leq e^{C'n+O(\log n)}$ for all $i$, Equation
\eqref{eq.partitions} and the above implies that for all $0 \leq k \leq
C n^2+O(n)$ we have:
\begin{eqnarray*}
|d_k| & \leq & \sum_{i=0}^k |a_i| c_{k-i} \\
& \leq &  e^{C' n +O(\log n)} \sum_{i=0}^k  p_{k-i} \\
& \leq &   e^{C' n+O(\log n)} k p_k \\
& \leq &  e^{C' n+O(\log n)} p_{C n^2+O(n)} \\
& \leq &  e^{C' n+O(\log n)}  e^{ \pi \sqrt{2C/3} n +O(\log n) } \\
& \leq & e^{C'' n+O(\log n)}
\end{eqnarray*}
where
$C''=C' + \pi \sqrt{2C/3}$. Since $||f||_1=\sum_{k=0}^{C n^2+O(n)} |d_k|$,
the result follows.
\end{proof}

\section{Some estimates}
\lbl{sec.hypgeom}

\subsection{The Lobachevsky function}
\lbl{sub.estimates}

In this largely independent section we will prove refined (and optimal)
estimates for the growth rate of the $R$-matrices. These estimates reveal
the close relationship between hyperbolic geometry and the asymptotics
of the quantum factorials. The main result in this section is a judicious
application of the Euler-MacLaurin summation formula.

As a warm-up, let us consider the {\em Lobachevsky function}
$$
\La(z)=-\int_0^z \log|2 \sin x| dx = \frac{1}{2} \sum_{n=1}^\infty 
\frac{\sin(2 n z)}{n^2}
$$
The Lobachevsky function is odd, with period $\pi$. 
Its graph for 
$z \in [0, \pi]$ is given by:
$$
\psdraw{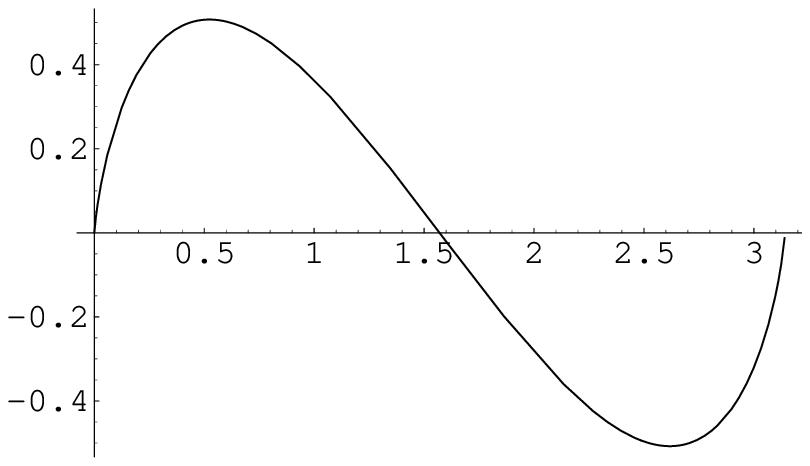}{2.5in}
$$
Notice also
that 
$v_8=8 \La(\pi/4) \approx 3.66386$ is the volume of the regular
ideal octahedron.

Recall that
if $f(q) \in \BZ[q^{\pm 1/2}]$, we denote by $\ev_n(f)$ the evaluation
of $f$ at $e^{2 \pi i/n}$.

Recall also the $q$-factorial $\{a \}!$ from Equation \eqref{eq.qfactorial}.
The next lemma discusses the asymptotics of the evaluation of
quantum factorials. The notation $O(\log n)$ below is a term which is bounded
by $C \log n$ for some constant $C$ {\em independent} of $\a$.

\begin{lemma}
\lbl{lem.1}
For every $\a \in (0,1)$ we have:
$$
\ev_n (\{ \lfloor \a n \rfloor \}!)=
\exp\left(-\frac{n}{\pi} \La(\pi \a) +O(\log n) \right) .
$$
\end{lemma}

\begin{remark}
\lbl{rem.subleading}
The proof reveals an asymptotic expansion of the form:
$$
\ev_n (\{ \lfloor \a n \rfloor \}!) \sim n^{\theta} 
\exp\left(-\frac{n}{\pi} \La(\pi \a) \right)\left(C_0 + \frac{C_1}{n} +
\frac{C_2}{C^2} + \dots\right)
$$
for explicitly computable constants $C_i$.
\end{remark}

\begin{proof}
Recall the {\em Euler-MacLaurin summation formula}, 
with error term (see for example, \cite[Chapt.8]{O}):
$$
\sum_{k=a}^b f(k)=\int_a^b f(x) dx + \frac{1}{2} f(a) + 
\frac{1}{2} f(b) + 
\sum_{k=1}^{m-1} \frac{B_{2k}}{(2k)!}( f^{(2k-1)}(b)-f^{(2k-1)}(a) )
+ R_m(a,b,f)
$$
where
$$
|R_m(a,b,f)| \leq (2-2^{1-2m}) \frac{|B_{2m}|}{(2m)!} 
\int_a^b |f^{(2m)}(x)| dx,
$$
and $B_k$ is the $k$th {\em Bernoulli number} given by the generating
series:
$$
\frac{x}{e^x-1}=\sum_{k=0}^\infty B_k \frac{x^k}{k!}.
$$
Applying the above formula to for $m=1$ to
$f(x)=\log|e^{2 \pi i x/n}-1|$, we have:
\begin{eqnarray*}
\log(\prod_{k=1}^{\lfloor \a n \rfloor } |e^{2 \pi i k/n}-1|) 
& = &  \frac{1}{2} (f(1)+f(\lfloor \a n \rfloor ))+
\int_0^{\lfloor \a n \rfloor } \log|e^{2 \pi i t/n}-1| dt +R_1(1,\a n,f) \\
& = &  \frac{1}{2} (f(1)+f( \a n )+
\int_0^{\a n } \log|e^{2 \pi i t/n}-1| dt +R_1(1,\a n,f) + \e(\a,n) \\
& = & \frac{1}{2} (f(1)+f(\a n))+
\frac{n}{\pi}  \int_{\pi/n}^{\pi \a} \log|2 \sin (  u)| u +R_1(1,\a n, f)
+\e(\a,n) 
\\
& = & \frac{1}{2} (f(1)+f(\a n))+
\frac{n}{\pi} \left(-\La(\pi \a) + \La(\frac{\pi}{n}) \right) + 
R_1(1,\a n, f) + \e(\a,n).
\end{eqnarray*}
Now, 
$$ 
\frac{1}{2} |f(1)+f(\a n)| \leq C \log n, \qquad \text{and}
\qquad | \e(\a,n) | \leq C''.
$$
Moreover, $f'(x)=-\frac{\pi}{n} \cot(\pi x/n)$ and $f''(x)=\frac{\pi^2}{n^2}
(\csc(\pi x/n))^2$, and $(\csc x)^2=1/x^2 + 1/3 + x^2/15 + O(x^3)$.
Thus $\int_{\pi/n}^{\a n} |\csc x| dx \leq C n^2$, and
$$
|R_1(1,\a n, f)| \leq C'. 
$$
Furthermore, using the asymptotic expansion of $\La(z)$ for $z \in (0, \pi)$:
$$
\La(z)= z -z \log(2 z) + \sum_{k=1}^\infty \frac{B_k}{2k} 
\frac{(2 z)^{2k+1}}{(2k+1)!},
$$
it follows that
$$
\frac{n}{\pi} |\La(\frac{\pi}{n})| \leq C'' \log n.
$$
The result follows.
\end{proof}

Consider the $R$-matrix  evaluated at $q=e^{2 \pi i/n}$,
$\ev_n (R(n;a,b,k))$. Recall from \eqref{eq.R+-} that $R$ is a ratio of 
$5$ quantum factorials.
Let us assume that $a=\lfloor \a n \rfloor $, $b=\lfloor \b n \rfloor $, 
$k=\lfloor \k n \rfloor $, where 
\begin{equation}
\lbl{eq.abk}
\a,\b,\k \in [0,1] \qquad 
0 \leq \b+\k \leq 1, \qquad 0 \leq \a-\k \leq 1.
\end{equation}

Let us define 
\begin{eqnarray*}
\lbl{eq.R}
r(n;\a,\b,\k) & :=& \ev_n(\log|R(n;\a n,\b n, \k n)|) \\ 
&=& \ev_n(\{b+k\}!) - \ev_n(\{b\}!) - \ev_n(\{k\}!)
+ \ev_n(\{a \}!) - \ev_n(\{a-k\}!) 
\end{eqnarray*}

Clearly,
$$
\max_{a,b,k} \ev_n(\log|R(n;a,b,k)|)=\max_{\a,\b,\k} r(n;\a,\b,\k).
$$
with the understanding that $\a,\b,\k$ satisfy \eqref{eq.abk}.

The next result gives the asymptotics of the $R$-matrix.

\begin{theorem}
\lbl{thm.R}
$$
\max_{\a,\b,\k}r(n;\a,\b,\k)=r(n;3/4,1/4,1/2)=
\frac{v_8}{2 \pi} n + O(\log n)
$$
where $v_8=8 \La(\pi/4)\approx 3.66386$ is the volume of the regular 
hyperbolic ideal octahedron.
\end{theorem}

Thus, asymptotically, the winning configuration is given by:
$$
\psdraw{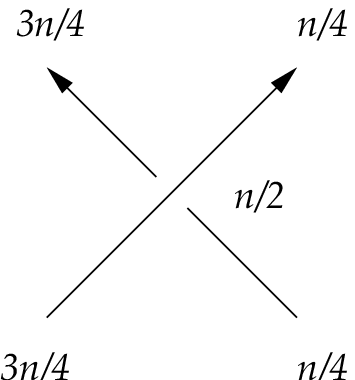}{1in}
$$

\begin{proof}
Lemma \ref{lem.1} and  thedefinition of $r$ imply that
\begin{eqnarray*}
r(n;\a,\b,\k) &=& \frac{n}{\pi} f(\a,\b,\k) + O(\log n)
\end{eqnarray*}
where
$$
f(\a,\b,\k)=-\La(\pi(\b+\k)) + \La(\pi\b)+\La(\pi\k)-\La(\pi\a)+\La(\pi(\a-\k))
$$
The domain of $f$ is a compact set, thus the maximum of $f$ exists.
Moreover, $f$ vanishes on the boundary, thus the maximum is one
of the critical points in the interior.
To find the critical points in the interior, let 
us set
$$
\za=e^{2 \pi i \a}, \qquad \zb=e^{2 \pi i \b}, \qquad
\zk=e^{2 \pi i \k}.
$$
Using the derivate of the Lobachevsky function
$$
\La'(x)=\log|e^{2 i x}-1|,
$$ 
it follows that
\begin{eqnarray*}
\frac{\pt f}{\pt \a} &=& \log|\za -1|-\log|\za \zk^{-1}-1|=0 \\
\frac{\pt f}{\pt \b} &=& \log|\zb \zk -1|-\log|\zb -1|=0 \\
\frac{\pt f}{\pt \k} &=& \log|\zb \zk -1| -\log|\zk -1|
-\log|\za \zk^{-1} -1|=0. 
\end{eqnarray*}
Thus, the critical points in the iterior are given by the solutions of:
\begin{eqnarray*}
|\za -1| &=& |\za \zk^{-1}-1| \\
|\zb \zk -1| &=& |\zb -1| \\
|\zb \zk -1| &=& |\zk -1| |\za \zk^{-1} -1|.
\end{eqnarray*}

Using Lemma \ref{lem.2circles} below, and the fact that $\za,\zb, \zk \neq 1$,
the first two equations imply that
$$
\za=\za^{-1} \zk \qquad \zb \zk =\zb^{-1}
$$
Thus, $\zk=\za^2=\zb^2$ and $\za=\pm \zb$. Plugging in the third equation
gives $(\za,\zb,\zk)=(-i,i,-1)$, i.e., $(\a,\b,\k)=(3/4,1/4,1/2)$.
Since
$$
\La(\pi/2)=0, \qquad \La(3 \pi/4)=-\La(\pi/4),
$$
it follows that
$$
f(3/4,1/4,1/2)=4 \La(\pi/4)=v_8/2.
$$
The result follows.
\end{proof}

\begin{remark}
\lbl{rem.octahedron}
The proof also reveals that $r(n;\a,\b,\k)=\frac{V(\a,\b,\k)}{2 \pi} n 
+ O(\log n)$, where $V(\a,\b,\k)$ is the volume of an ideal octahedron
with vertices $0,1,\infty,\za,\zb,\zk$.
\end{remark}

\begin{lemma}
\lbl{lem.2circles}
If $z,w$ are complex numbers that satisfy $|z|=|w|=1$ and $|1-z|=|1-w|$,
then $z=w^{\pm}$.
\end{lemma}

\begin{proof}
Let us define 
$$
C_{u_0,r}:=\{ u \in \BC \,\, | \,\, |u-u_0|=r>0 \}.
$$
Then $C_{u_0,r}$ is a circle with center $u_0$ and radius $r$. Fixing $w$, it 
follows that $z \in C_{0,1} \cap C_{1,|1-w|}$. The intersection of two circles
is two points, and since $w$ and $w^{-1}=\bar{w}$ both lie in the intersection,
the result follows.
\end{proof}

\begin{remark}
\lbl{rem.othercases}
The same bound in Theorem \ref{thm.R} holds for the $R_-$ matrix of
Equation \eqref{eq.R+-}.
\end{remark}

\subsection{Proof of Theorem \ref{thm.4}}
\lbl{sub.thm4}

Recall that the colored Jones function is given by the state-sum of
Equation \eqref{eq.statesum}, where the summand $F(n,\bk)$ is a product
of local $R$-matrices, one for each crossing of $K$. 
Theorem \ref{thm.R} implies that
$$
|\ev_n(F(n,\bk))| \leq e^{ v_8/(2\pi) n c + O(\log n)}
$$
for each $\bk$, where the error term is bounded independent of $\bk$.
Since $\bk$ takes $O(n^4)$ values, Theorem \ref{thm.4} follows.

\section{The $q$-holonomic point of view}
\lbl{sec.qholonomic}

\subsection{Bounds on $l^1$-norm and Mahler measure of $q$-holonomic functions}
\lbl{sub.qholo}

The main result of \cite{GL} is that for every knot $K$, the functions $J_K$
and $C_K$ are $q$-holonomic.
Recall that a sequence $f:\BN\longto\BQ(q)$ is $q$-{\em holonomic} if 
satisfies a $q$-{linear difference equation}. In other words,
there exists a natural number $d$ and rational functions 
$a_j(u,v) \in \BQ(u,v)$ for $j=0,\dots,d$ with $a_d \neq 0$ such that
for all $n \in \BN$ we have:
$$
\sum_{j=0}^d a_j(q^n,q) f(n+j)=0.
$$

In this section we observe that $q$-holonomic functions satisfy {\em a priori}
upper bounds on their degrees and (under an integrality assumption)
on their $l^1$-norm. As a simple corollary, we obtain
an independent proof of Theorems \ref{thm.1} and \ref{thm.3}.

\begin{definition}
\lbl{eq.integral}
We say that a sequence $f:\BN\longto\BZ[q^{\pm}]$ is $q$-{\em integral
holonomic} if it
satisfies an {\em integral
$q$-difference equation} as above with
$a_d=1$ and $a_j(u,v) \in \BZ[u,v]$.
\end{definition}

Although $J_K$ takes values in $\BZ[q^{\pm}]$, it is known
for example that $J_{4_1}$ is not $q$-integral holonomic. On the other hand,
it is known that $C_K$ is $q$-integral holonomic for all twist knots;
see \cite{GS}.  

\begin{question}
\lbl{que.1}
Is it true that $C_K$ is $q$-integral holonomic for every knot $K$?
\end{question}

For a Laurent polynomial $f(q)=\sum_{k=m}^M a_k q^k$, with $a_m a_M \neq 0$,
let us define $\text{degmax}_q(f)=M$ and $\text{degmin}_q(f)=m$.

\begin{theorem}
\lbl{thm.L1bound}
$\mathrm{(a)}$
If $f:\BN\longto \BZ[q^{\pm}]$ is $q$-holonomic, then for all $n$ we have:
$$
\text{maxdeg}_q(f(n))=O(n^2)   \qquad
\text{and} \qquad \text{mindeg}_q(f(n)) =O(n^2).
$$
\newline
$\mathrm{(a)}$ If $f$ is $q$-integral holonomic, then for all $n$ we have:
$$
||f(n)||_1 \leq C^n
$$
for some constant $C$. In particular,
$$
M(f(n)) \leq C^n
$$
and
$$
\limsup_{n \to \infty} \frac{\log|\ev_{\a,n}(f(n))|}{n}
\leq C
$$
for all $\a \in \BR$.
\end{theorem}

\begin{proof}
For the first claim in (a),
let us assume without loss of generality that
$$
a_d(q^n,q) f(n+d) = -a_{d-1}(q^n,q) f(n+d-1) - \dots a_0(q^n,q) f(n)
$$
where $a_j(Q,q) \BQ[Q,q]$ are polynomials in $Q,q$. Choose $C'$ so that
\begin{itemize}
\item
$\text{maxdeg}_q a_j(q^n,q) \leq 2 C'(n+d)$ for all $j=0,\dots,d-1$, and
\item
$\text{maxdeg}_q f(n) \leq C'(n+1)^2$ for $n=0,\dots,d-1$.
\end{itemize}
We will prove by induction on $n$ that $\text{maxdeg}_q f(n) \leq C'(n+1)^2$.
By assumption, it is true for $n=0, \dots, d-1$.
Then, by induction we have:
\begin{eqnarray*}
\text{maxdeg}_q f(n+d) & \leq & 
\text{maxdeg}_q a_d(q^n,q) + \text{maxdeg}_q f(n+d) \\
&=& \text{maxdeg}_q (a_d(q^n,q) f(n+d)) \\
& \leq & \max_{0 \leq j < d} \text{maxdeg} (a_j(q^n,q) f(n+j)) \\
&=& \max_{0 \leq j < d} 2 C' (n+j+1) +C'(n+j+1)^2 \\
&=&  2 C'(n+d) +  C'(n+d)^2 \\
& < & C'(n+d+1)^2.
\end{eqnarray*}
The second claim in (a) follows similarly.

For (b), let $c_j=||a_j(Q,q)||_1$ for $j=0, \dots, d-1$, and choose $C$
so that 
\begin{itemize}
\item
$C^d \leq c_{d-1} C^{d-1} + \dots + c_0 C^0$, and
\item
$||f(n)||_1 \leq C^n$ for $n=0, \dots, d-1$.
\end{itemize}
Then, it is easy to see by induction that (b) holds for all $n$.
\end{proof}

As advertised above, Theorem \ref{thm.L1bound} gives an alternative proof
of Theorem \ref{thm.1} and \ref{thm.3}, under the assumption that Question
\ref{que.1} has a positive answer.
However, the explicit upper bounds in terms of the number
of crossings cannot be obtained from Theorem \ref{thm.L1bound},
unless we know something more about the $q$-difference equation of the
colored Jones function. Moreover, Theorem \ref{thm.2} cannot be obtained
from general theory of asymptotics of solutions of $q$-difference equations,
since a typical solution would be growing exponentially, 
even for small positive
$\a$. With additional assumptions on the shape of the $q$-difference
equations, the first author can show Theorem \ref{thm.2}. The proof
was replaced with the one of the present paper.

\subsection{Bounds for higher rank groups}
\lbl{sub.high}

In \cite{GL}, we considered the colored Jones function 
$$
J_{\fg,K}: \La_w \longto \BZ[q^{\pm}]
$$ 
of a knot $K$, where $\fg$ is a {\em simple Lie algebra} with {\em 
weight lattice} 
$\La_w$. In the above refernce, the authors proved that $J_{\fg,K}$
is a $q$-holonomic function, at least when $\fg$ is not $G_2$.
For $\fg=\mathfrak{sl}_2$, $J_{\mathfrak{sl}_2,K}$ is the colored Jones
function $J_K$ discussed earlier.

In \cite{GL} , the authors gave  state-sum formulas for $J_{\fg,K}$
similar to \eqref{eq.statesum} where the summand takes values in 
$\BZ[q^{\pm 1/D}]$, where $D$ is the size of the center of $\fg$.

The methods of the present paper give an upper bound for
the growth-rate of the $\fg$-colored Jones function. More precisely,
we have:

\begin{theorem}
\lbl{thm.gJones}
For every simple Lie algebra $\fg$ (other than $G_2$)  there exists a 
constant $C_{\fg}$ such that for every knot with $c+2$ crossings, and
every $\a >0$ and every $\l \in \La_w$, we have:
$$
\limsup_{n \to \infty} 
\frac{\log|\ev_{\a,n}(J_{\fg,K}(n \l))|}{n}
\leq C_{\fg} c .
$$
\end{theorem}

The details of the above theorem will be explained in a subsequent
publication.

\subsection{Acknowledgement} 

The authors wish to thank I. Agol, D. Boyd, N. Dunfield, D.
Thurston and D. Zeilberger for many enlightening conversations.

\appendix

\section{The volume conjecture for the Borromean rings}
\lbl{sec.borromean}

It is well-known that the complement of the Borromean rings $B$ can be
geometrically identified by gluing two regular ideal octahedra;
\cite{Th}. As a result, the volume $V(B)$ of $B$
is given by $2 v_8$.

If $L$ is a link with a distinguished component (to be broken), then
one may define the colored Jones function $J_L(n)$ to be the invariant
of the $(1,1)$-tangle obtained by breaking the distinguished component
of $L$ and coloring all components of the tangle by the $n$-dimensional
irreducible representation of $\mathfrak{sl}_2$. In general, $J_L(n)$
depends on the link and its distinguished component. In the case of the
Borromean rings $B$ though, due to symmetry, we may choose any component
as the distinguished one. Habiro uses the notation $\tilde{J}_L(n)$ 
for $J_L(n)$.

The next theorem confirms the volume conjecture for the Borromean rings.

\begin{theorem}
\lbl{thm.borromean}
If $J_B(n)$ denotes the colored Jones function of the Borromean rings
$B$, then
\begin{equation}
\lbl{eq.borro}
\lim_{n \to \infty} \frac{\log|\ev_n (J_B(n))|}{n}
=\frac{1}{2 \pi} V(B).
\end{equation}
\end{theorem}

\begin{proof}
Using  Habiro's formula for $\tilde J_L$ of the Borromean ring
\cite{H2},  one has
$$
J_B(n)= \sum_{l=0}^{N-1}(-1)^l\, \frac{\{n\}^2 \left( \prod_{j=1}^{l} \{n+j\}
\{n-j\}\right)^3}{\left(\prod_{j=l+1}^{2l+1}\{j\}\right)^2}.
$$
When $v= e^{i \pi/n}$, on has $\{j\}= 2i \sin \frac{j\pi}{n}$, which
is 0 exactly when $j$ is divisible by $n$. Hence if $2l+1 <n$,
then the denominator of the term in the above sum is never 0,
while the numerator is 0, since it has 2 factors $\{n\}$. On the
other hand, if $2l+1 >n$, then the denominator has 2 factors
$\{n\}$, which would cancel with the 2 same  factors of the
numerator. Hence when evaluating at $v= e^{i \pi/n}$ one can assume
that $2l +1 \ge n$, or $l> n/2-1$:
$$ 
\ev_n(J_B(n))= \sum_{n>l>n/2-1}(-1)^l\, \frac{ \left( \prod_{j=1}^{l} \{n+j\}
\{n-j\}\right)^3}{\left(\prod_{j=l+1}^{n-1}\{j\}\,
\prod_{j=n+1}^{2l+1}\{j\}\right)^2}.
$$
Note that when $v= e^{i \pi/n}$, one has $\{n+j\}=-\{j\} = -2i\sin
(j\pi/n)$. Hence a simple calculation shows that
$$ 
\ev_n(J_B(n))= \sum_{n>l>n/2-1} 
\frac{ \left( \tau_{1,l}\right)^3}{\tau_{l+1,n-1}\, \tau_{n+1,2l+1}},
$$
where
$$
\tau_{p,l} := \prod_{j=p}^{l} 4 \sin^2 (j\pi/n).
$$
The following properties of $\tau$ are easy to verify

\begin{lemma}
\lbl{lem.tau} 
We have:
\begin{eqnarray}
\tau_{1,m} & = & \tau_{n-m,n-1} \qquad \text{for} \qquad
0<m<n
\lbl{a} \\
\tau_{n+1,m} & = & \tau_{2nn-m,n-1} \qquad \text{for} \qquad
n<m<2n
\lbl{b} \\
\tau_{1,n-1} &=& n^2
\lbl{c} \\
\tau_{1,m} \,\tau_{1,n-m-1} &=& n^2 \qquad \text{for} \qquad
0<m<n.
\lbl{d}
\end{eqnarray}
\end{lemma}

From the Lemma it follows that $\tau_{l+1,n-1}= n^2/\tau_{1,l}$ and
$\tau_{n+1,2l+1}= \tau_{2n-l-1,n-1} = n^2/\tau_{1,2n-l-2}$. 
Hence
$$ 
\ev_n(J_B(n))= n^4 \sum_{n>l>n/2-1} 
\left( \tau_{1,l}\right)^4 \, \tau_{1,2n-2l-2}.
$$
Let $k= 2l+1-n$. Then $n>k\ge 0$, $k+n \equiv 1 \pmod{2}$, and $l=
(n+k-1)/2$. From the above Lemma one has that
$$
\tau_{1, 2n-2l-2} = \frac{n^2}{\tau_{1,k}},
$$
$$
\tau_{1,l} = \frac{n^2}{\tau_{n-l-1}}.
$$
Thus,
$$ 
(\tau_{1,l})^4 = (\tau_{1,l})^2 \, \left(\frac{n^2}{\tau_{n-l,-1}}\right)^2,
$$
and

\begin{equation} 
\lbl{st}
\ev_n(J_B(n))= n^{10}\sum_{n>k\ge 0, k+n\equiv
1\pmod2}\frac{(\tau_{1, (n+k-1)/2})^2}{\tau_{1,(n-k-1)/2}
\tau_{1,k}}.
\end{equation}
Using the notation of Section \ref{sec.hypgeom}, we have:
$$
\frac{(\tau_{1, (n+k-1)/2})^2}{\tau_{1,(n-k-1)/2}
\tau_{1,k}}=|\ev_n(R(n;a,b,k))|^2
$$
where $a= \lfloor (n+k-1)/2 \rfloor$ and $b= \lfloor (n-k-1)/2 \rfloor$. 
There are
less than $n$ terms, hence using Theorem \ref{thm.R} one has that
$$ 
\ev_n(J_B(n)) < n^{11} e^{2 r(n;3/4,1/4,1/2)}.
$$
Therefore the limsup of the left hand side of \eqref{eq.borro} is 
$2 v_8/(2 \pi)$.

On the other hand, all the terms in the sum of (\ref{st}) are positive, 
hence the sum is at least the biggest term, which is when 
$k=\lfloor n/2 \rfloor$. 
The limit when we retain only this term
is easily seen to be equal to $(2 v_8)/(2 \pi)$, 
which is the same as the limsup.
\end{proof}

\section{The volume conjecture for torus knots}
\lbl{sec.torusknots}

Let $T(a,b)$ denote the $(a,b)$-torus knot for
co-prime integers $a,b$ ($a,b >1$). Although $T(a,b)$ is not a
hyperbolic knot, the {\em volume function} $V(l,m)$ on the deformation
variety can be defined; see \cite{CCGLS}.  Let us discuss this first.

\begin{lemma}
The restriction of the volume form on the deformation variety,
i.e., the zero set of the A-polynomial is equal to 0.
\end{lemma}
\begin{proof} 
For a torus knots the A-polynomial is either of the
form $l\pm m^c$ for some integer $c$, or the product of 2 factors
of that forms. It is easy to check that the restriction of the
volume form on any such factor is equal to 0.
\end{proof}

Thus the volume function $V(l,m)$ on the deformation variety must
be constant. Because the Gromov norm of the knot complement is 0,
one should define $V(l,m)=0$ for every $l,m$. Then the generalized
volume conjecture can be proved easily in this case:

\begin{proposition}
\lbl{prop.torus}
For the torus knot $K=T(a,b)$ and for any real number $\a$, one has
\begin{equation}
 \lim_{n \to \infty} \frac{\log|\ev_{\a,n}(J_K(n))|}{n}= 0.
\end{equation}
\end{proposition}

\begin{proof}
There are 2 cases: $\a$ is an integer, or $\a$ is not. The first
case is actually much more difficult, this is the usual volume
conjecture and it has been proved by Kashaev and Tirkkonen in 
\cite{KT}. Let us consider the easier case, when $\a$ is not
an integer.

Note that $v= \exp(\pi i \a/n)$, and if $\a$ is not an integer one
has  $v^n-v^{-n}= 2 i \sin \pi \a/n \neq 0$.

The colored Jones polynomial was calculated by Morton in \cite{Mo}:

\begin{equation}
\lbl{equ.torus}
J_{T(a,b)(n)} = \frac{v^{-ab(n^2-1)/2}}{v^n -v^{-n}} \sum_{k= 1-n,
k+n\equiv 1 \pmod2}^{n-1}(v^{2abk^2+ 2ka +2kb +1} - v^{2abk^2+ 2ka
-2kb -1}). 
\end{equation}

The sum contains $n+1$ terms, each by absolute value is less than
or equal to 2. Hence
$$ 
|\ev_{\a,n}(J_K(n))| 
< \frac{2(n+1)}{\sin (\pi \a/n)}.
$$ 
Thus, the limsup is less than or equal to 0. 
The argument in Section \ref{sec.thm2}
shows that the $\liminf$ is greater than or equal to 0. The result follows.
\end{proof}

\begin{remark}
\lbl{rem.mur}
In \cite{M2}, H. Murakami discusses the Generalized Voume Conjecture for
the torus knots and angles $\a$ with nonzero imaginary part.
\end{remark}

\ifx\undefined\bysame
    \newcommand{\bysame}{\leavevmode\hbox
to3em{\hrulefill}\,}
\fi

\end{document}